\documentclass[a4paper,oneside,11pt]{article}

\usepackage{amsmath}
\usepackage{amssymb}
\usepackage{amsfonts}
\usepackage{theorem}
\usepackage{bbm}

\usepackage[T1]{fontenc}
\usepackage{times}
\usepackage[cp1250]{inputenc}

\newtheorem{dfn}{Definition}[section]
\newtheorem{tw}[dfn]{Theorem}
\newtheorem{prop}[dfn]{Proposition}
\newtheorem{rem}[dfn]{Remark}

\newtheorem{lem}[dfn]{Lemma}

\usepackage{anysize}
\usepackage{verbatim}
\usepackage{array}
\usepackage{enumerate}
\usepackage{amssymb} 
\usepackage{amsmath}
\usepackage{fancyhdr}
\usepackage{euscript}
\usepackage{latexsym}
\usepackage{graphicx}

\author{Micha\l \ Barski  \\ \small  Faculty of Mathematics, Cardinal Stefan Wyszy\'nski University in Warsaw, Poland\\
\small Faculty of Mathematics and Computer Science, University of Leipzig, Germany\\ \small{\it Michal.Barski@math.uni-leipzig.de} \bigskip \\
\\
Jerzy Zabczyk
\\ \small Institute of Mathematics, Polish Academy of
Sciences,
     Warsaw,  Poland\\ \small{\it zabczyk@impan.pl}}

\title{\bf Completeness of bond market driven by L\'evy process\thanks{Research supported by Polish KBN Grant
P03A 034 29 ,,Stochastic evolution equations driven by L\'evy
noise''} }

\begin{document}
\baselineskip=1.1\baselineskip \maketitle

\begin{abstract}
The completeness problem of the bond market model with the random
factors determined by a Wiener process and Poisson random measure is
studied. Hedging portfolios use bonds with maturities in a
countable, dense subset of a finite time interval. It is shown that
under natural assumptions the market is not complete unless the
support of the L\'evy  measure consists of a finite number of
points. Explicit constructions of contingent claims which can not be
replicated are provided.
\end{abstract}

\noindent
\begin{quote}


\end{quote}

\bigskip
\section{Introduction}
\quad\ Tradeable bonds are specified by a set of their maturities,
which potentially can be composed of infinitely many points. Thus
the bond market consists of infinitely many assets and this is a
significant difference with respect to classical market of a finite
number of stocks. This is a reason why the bond market models are
not covered by the classical theory and economic problems, like
completeness, have to be studied anew.

The problem of bond market completeness was treated in many
different contexts depending on the model settings as well as on the
definition of completeness. A classical question of the market
completeness is to judge if it is possible to replicate any bounded
random variable $X$, i.e. to find a portfolio which is equal to $X$
at the final time. However, it is sometimes difficult to solve this
problem in the set of all bounded random variables and thus another
spaces are also considered, for example $L^2(\Omega)$ or even more
exotic ones. In Taflin \cite{Taf} it is shown that the model driven
by the infinite dimensional Wiener process is not complete in the
class $D_{0}:=\bigcap_{p>1}L^p(\Omega)$. In Carmona, Tehranchi
\cite{Car} it is shown that each random variable which is of a
special form can be replicated.

Another question connected with the notion of completeness is that
of existence of a unique martingale measure. Contrary to the finite
dimensional market this property of the model, in general, is not
equivalent to completeness. As it was shown in Bj\"ork et. al
\cite{BjoKabRunMas} and \cite{BjoKabRun} in a jump diffusion model
uniqueness of the martingale measure is equivalent to the
approximate completeness, i.e. for any random variable $X\in
L^2(\Omega)$ there exists a sequence of random variables $\{X_n\}$
which converges to $X$ in $L^2(\Omega)$ s.t. each element of the
sequence can be replicated.

It was shown in Baran, Jakubowski, Zabczyk \cite{BarJakZab} that a
model driven by the infinite dimensional Wiener process is not
complete, i.e. there exists a bounded random variable which can not
be replicated. In this paper we focus on a finite dimensional noise
with jumps and for simplicity assume that it is given by the one
dimensional Wiener process and Poisson random measure. We consider
model with a finite time interval $[0,T^\ast]$. Each bond is
specified by its maturity $T$ and usually it is assumed that
maturity can by any number from $[0,T^\ast]$. We adopt the setting
of Eberlein, Jacod, Raible \cite{EbJacRaib} and consider bonds with
maturities in a dense, countable subset of $[0,T^\ast]$ denoted by
$J$. This set consists of all bonds' maturities which can be
involved in the portfolios construction. A bond with maturity $T$
and the price process $P(\cdot,T)$ can be used by a trader if and
only if $T\in J$. The completeness problem with the use of bonds
with maturities in $J$ can be formulated in two ways:
\begin{enumerate}[1)]
\item Does there exist a unique equivalent measure $Q$ such that the discounted
prices of bonds $\hat{P}(\cdot,T)$ are $Q$-local martingales for
each $T\in J$?
\item Can arbitrary $\mathcal{F}_{T^\ast}$- measurable random variable,
satisfying some regularity assumptions, be replicated with the use
of bonds with maturities in $J$?
\end{enumerate}
Analogous formulations to $(1)$ and $(2)$ for finite number of
stocks are equivalent - at least for a wide class of stock market
models. However, as it was shown in \cite{BjoKabRunMas} and
\cite{BjoKabRun} they can no longer be equivalent if we examine bond
market with infinite number of assets. The problem of completeness
with the use of bonds with maturities in $J$ was originally
formulated in \cite{EbJacRaib}, where it was treated in the sense of
the formulation $(1)$. It was shown that under appropriate
assumptions there exists exactly one martingale measure. In this
paper we study the problem of completeness in the sense of the
formulation $(2)$. This approach requires a precise definition of
portfolios which can be used by traders, see Section \ref{sec Bond
market model}. We identify prices of bonds with elements of a Banach
space $B$ consisting of all bounded sequences with the supremum
norm. The trader's position is identified with an element of $l^1$ -
a subspace of the dual space $B^\ast$. The self-financing condition
is expressed by the fact that portfolio's value is an integral of
the $l^1$-valued strategy with respect to the bond price process.

The general idea in the solution of the completeness problem is to
examine the possibility of representing any martingale as a certain
stochastic integral with $l^1$-valued integrand. The key tools used
for this purpose are the representation theorem for local
martingales, which comes from Kunita \cite{Kun}, and a version of
theorem solving the so called problem of moments. The last one
provides necessary and sufficient conditions for the existence of a
linear, bounded functional satisfying certain conditions. Generally
speaking we apply this theorem to the real and vector-valued
functions defined on the support of the L\'evy measure. Our main
result states that every market model with the L\'evy measure having
a concentration point is incomplete. We provide an explicit
construction of a bounded random variable which can not be
replicated. If there is no concentration point we prove
incompleteness under additional assumptions in the class of square
integrable or bounded random variables. In the case when the L\'evy
measure has a finite support and the model satisfies additional
assumptions we prove completeness in the class of integrable random
variables. This result is similar to Theorem 5.6 in \cite{BjoKabRun}
but requires weaker assumptions.

The paper is organized as follows: in Section \ref{sec Local
martingales representation} we recall basic facts on stochastic
integrals and formulate the representation theorem for local
martingales; Section \ref{sec Bond market model} contains a
description of the model and definition of portfolios; in Section
\ref{sec Completeness} we present the main results - this section is
divided into three parts with respect to the properties of the
L\'evy measure.

\section{Local martingales representation}\label{sec Local martingales representation}

We will consider a c\`adl\`ag version of the L\'evy process
$Z=\{Z(t); t\in[0,T^\ast]\}$ defined on a probability space
$(\Omega,\mathcal{F},P)$. Let $N$ be the associated jump measure
\begin{gather*}
N(t,A):=\sharp\{s\in[0,t]: \triangle Z(s)=Z(s)-Z(s-)\in A\}, \quad
t\in[0,T^\ast], \ A\subseteq\mathbb{R}.
\end{gather*}
If $A$ is such that $0\notin\bar{A}$ then $N(t,A)$ is integrable and
its expectation can be written in the form
\begin{gather*}
\mathbf{E}N(t,A)=t \ \nu(A), \quad t\in[0,T^\ast], \ 0\notin\bar{A}.
\end{gather*}
The measure $\nu$ above, called the L\'evy measure, is such that
\begin{gather*}
\int_{\mathbb{R}}\mid x\mid^2\wedge \ 1 \ \nu(dx)<\infty, \quad
\nu(\{0\})=0.
\end{gather*}
The compensated jump measure $\tilde{N}$ is defined by
\begin{gather*}
\tilde{N}(t,A):=N(t,A)-t\nu(A)\quad t\in[0,T^\ast], 0\notin\bar{A}.
\end{gather*}
It is known that $Z$ can be decomposed into the following
L\'evy-It\^o form, see \cite{Apl}, Th. 2.4.16,
\begin{gather*}
Z(t)=at+bW(t)+\int_{0}^{t}\int_{\mid x
\mid<1}x\tilde{N}(ds,dx)+\int_{0}^{t}\int_{\mid x \mid\geq1}x
N(ds,dx);\qquad t\in[0,T^\ast],
\end{gather*}
where $a\in\mathbb{R}$, $b\geq 0$ and $W$ is a standard Wiener
process adapted to the filtration
\begin{gather*}
\mathcal{F}_t:=\sigma\{Z(s): 0\leq s\leq t\},\qquad t\in[0,T^\ast],
\end{gather*}
generated by $Z$.

\medskip
\noindent In order to formulate the representation theorem below, we
briefly present description of the class of integrable processes
with respect to $W$ and $\tilde{N}$. We follow notation used in
\cite{Kun}.

\medskip
\noindent The process $\phi=(\phi(\omega,t))$ is integrable with
respect to the Wiener process if it is predictable and satisfies
integrability condition
\begin{gather*}
\int_{0}^{T^\ast}\mid\phi(s)\mid^2 ds<\infty, \qquad \ P-\text{a.s.}
.
\end{gather*}
This class of processes is denoted by $\Phi$. For any $\phi\in\Phi$
the integral
\begin{gather*}
\int_{0}^{t}\phi(s)dW(s):=\int_{0}^{T^\ast}\phi(s)\mathbf{1}_{[0,t]}(s)dW(s)
\end{gather*}
is well defined and the process $\int_{0}^{\cdot}\phi(s)dW(s)$ is a
continuous locally square integrable martingale.\\

\noindent The process $\psi=(\psi(\omega,s,x))$ is called
predictable if it is $\mathcal{P}\otimes\mathcal{B}(\mathbb{R})$
measurable, where $\mathcal{P}$ is a predictable sigma-field. If
$\psi$ satisfies condition
\begin{gather}\label{war na Psi1}
\int_{0}^{T^\ast}\int_{\mathbb{R}}\mid\psi(s,x)\mid\nu(dx)ds<\infty,
\qquad P-\text{a.s.},
\end{gather}
then the integral
\begin{gather*}
\int_{0}^{T^\ast}\psi(s,x)\tilde{N}(dx,ds)=\int_{0}^{T^\ast}\psi(s,x)N(ds,dx)-\int_{0}^{T^\ast}\psi(s,x)\nu(dx)ds
\end{gather*}
is well defined and the process
$\int_{0}^{\cdot}\psi(s,x)\tilde{N}(ds,dx)=\int_{0}^{T^\ast}\psi(s,x)\mathbf{1}_{(0,\cdot]}(s)\tilde{N}(ds,dx)$
is a local martingale. The class of predictable processes satisfying
\eqref{war na Psi1} is denoted by $\Psi_1$. \\

\noindent If a predictable process $\psi$ satisfies condition
\begin{gather}\label{war na Psi2}
\int_{0}^{T^\ast}\int_{\mathbb{R}}\mid\psi(s,x)\mid^2\nu(dx)ds<\infty
\qquad P-\text{a.s.}
\end{gather}
then the integral $\int_{0}^{T^\ast}\psi(s,x)\tilde{N}(ds,dx)$ is
constructed with the use of simple processes which converge to
$\psi$ in $L^2$. In this case
$\int_{0}^{\cdot}\psi(s,x)\tilde{N}(ds,dx)=\int_{0}^{T^\ast}\psi(s,x)\mathbf{1}_{(0,\cdot]}(s)\tilde{N}(ds,dx)$
is a locally square integrable martingale. A class of predictable
processes satisfying \eqref{war na Psi2} is denoted by $\Psi_2$.
\\
\noindent A class of all predictable processes which satisfy
conditions
\begin{gather*}
\psi\mathbf{1}_{\{\mid\psi\mid>1\}}\in\Psi_1 \quad\text{and}\quad
\psi\mathbf{1}_{\{\mid\psi\mid\leq1\}}\in\Psi_2
\end{gather*}
will be denoted by $\Psi_{1,2}$. In other words $\psi\in\Psi_{1,2}$
if and only if
\begin{gather*}
\int_{0}^{T^\ast}\int_{\mathbb{R}}\mid\psi(s,x)\mid^2\wedge\mid\psi(s,x)\mid
\nu(dx)ds<\infty.
\end{gather*}

\noindent For any $\psi\in\Psi_{1,2}$ the integral
\begin{align*}
\int_{0}^{T^\ast}\psi(s,x)\tilde{N}(ds,dx)&=\int_{0}^{T^\ast}\psi(s,x)\mathbf{1}_{\{\mid\psi(s,x)\mid>1\}}(s,x)\tilde{N}(ds,dx)\\[2ex]
&+\int_{0}^{T^\ast}\psi(s,x)\mathbf{1}_{\{\mid\psi(s,x)\mid\leq1\}}(s,x)\tilde{N}(ds,dx)
\end{align*}
is well defined and it is a local martingale as a function of the
upper integration limit.\\

\noindent The next theorem comes from \cite{Kun}.
\begin{tw}\label{tw Kunity}
Let $M$ be an $(\mathcal{F}_t)$-local martingale. Then there exist
$\phi\in\Phi$ and $\psi\in\Psi_{1,2}$ satisfying
\begin{gather}\label{reprezentacja Kunity}
M_t=M_0+\int_{0}^{t}\phi(s)dW(s)+\int_{0}^{t}\int_{\mathbb{R}}\psi(s,x)\tilde{N}(dx,ds).
\end{gather}
Moreover, the pair $(\phi,\psi)$ is unique i.e., if
$(\phi^{'},\psi^{'})$ satisfies \eqref{reprezentacja Kunity} then
\begin{gather*}
\phi=\phi^{'} \ w.r.t. \ \ P\otimes\lambda- \ \text{a.s.}  \quad
\text{and} \quad \psi=\psi^{'}   \ w.r.t. \ \
P\otimes\lambda\otimes\nu- \ \text{a.s.} ,
\end{gather*}
where $\lambda$ is the Lebesgue measure on $[0,T^{\ast}]$.
\end{tw}

\section{Bond market model}\label{sec Bond market model}

We begin description of the model by specifying the dynamics of the
forward rate
\begin{gather}\label{rownanie na f}
df(t,T)=\alpha(t,T)dt+\sigma(t,T)dW(t)+\int_{\mathbb{R}}\gamma(t,x,T)
N(dt,dx),\quad  t,T\in[0,T^\ast].
\end{gather}
The coefficients are assumed to be predictable and satisfy the
following integrability conditions
\begin{gather*}
\int_{0}^{T^\ast}\int_{0}^{T^\ast}\mid \alpha(t,T)\mid dTdt<\infty,
\ \int_{0}^{T^\ast}\int_{0}^{T^\ast}\mid \sigma(t,T)\mid^2
dTdt<\infty, \\[2ex]
\int_{0}^{T^\ast}\int_{0}^{T^\ast}\int_{\mathbb{R}}\mid
\gamma(t,x,T)\mid \nu(dx)dTdt<\infty,
\end{gather*}
where all the inequalities above hold $P$-a.s.. We put
\begin{gather}\label{war na zerowanie wspolczynnikow}
\alpha(t,T)=0,\quad\sigma(t,T)=0,\quad\gamma(t,x,T)=0 \quad
\text{for} \quad t>T, \quad \forall x\in\mathbb{R}.
\end{gather}

\noindent
 The value at time $t$ of a bond paying $1$ at maturity
$T\in[0,T^\ast]$ is defined by
\begin{gather}\label{rownanie def na P}
P(t,T):=e^{-\int_{t}^{T}f(t,s)ds},\qquad t,T\in[0,T^\ast].
\end{gather}
 \noindent The evolution of the money in the savings account is
given by
\begin{gather*}
dB(t)=r(t)B(t)dt,\qquad t\in[0,T^\ast],
\end{gather*}
where $r(t):=f(t,t)$ is the short rate. In virtue of \eqref{war na
zerowanie wspolczynnikow} we have equality $f(t,T)=f(T,T)$ for
$t>T$. Indeed, for $t\in[T,T^\ast]$, we have
\begin{align*}
f(t,T)&=f(0,T)+\int_{0}^{t}\alpha(s,T)ds+\int_{0}^{t}\sigma(t,T)dW(s)\\[2ex]
&=f(0,T)+\int_{0}^{T}\alpha(s,T)ds+\int_{0}^{T}\sigma(t,T)dW(s)\\[2ex]
&=f(T,T).
\end{align*}
This relation implies the following equality
\begin{align*}
P(t,T)&=e^{-\int_{t}^{T}f(t,s)ds}=e^{-\int_{t}^{T}f(s,s)ds}=e^{-\int_{t}^{T}r(s)ds}\\[2ex]
&=e^{\int_{T}^{t}r(s)ds}=P(T,T)e^{\int_{T}^{t}r(s)ds}, \qquad
\text{for} \ \ t\in[T,T^\ast],
\end{align*}
which corresponds to the fact that the holder of a bond transfers
his money automatically to the bank account after the bond's
expiration date.

\noindent The discounted value of the bond
$\hat{P}(t,T):=B(t)^{-1}P(t,T)$ with maturity $T$ is thus given by
\begin{gather*}
\hat{P}(t,T)=P(t,T)\
e^{-\int_{0}^{t}r(s)ds}=e^{-\int_{t}^{T}f(t,s)ds}\cdot
e^{-\int_{0}^{t}f(t,s)ds}=e^{-\int_{0}^{T}f(t,s)ds}, \quad
t,T\in[0,T^\ast].
\end{gather*}
As a consequence, the discounted value of the bond
\begin{gather*}
\hat{P}(t,T)=e^{-\int_{0}^{T}f(t,s)ds}=e^{-\int_{0}^{T}f(s,s)ds}=e^{-\int_{0}^{T}r(s)ds},
\qquad t\in[T,T^\ast],
\end{gather*}
is constant after its expiration date.

\noindent Putting
\begin{align*}
A(t,T):=-\int_{t}^{T}\alpha(t,s)ds\\
S(t,T):=-\int_{t}^{T}\sigma(t,s)ds\\
G(t,x,T):=-\int_{t}^{T}\gamma(t,x,s)ds\\
\end{align*}
one can check that $P$ satisfies the following equation (see
Proposition 2.2. in \cite{BjoKabRun}):
\begin{align}\label{rownanie na P}\nonumber
dP(t,T)=P(t-,T)\bigg(\Big( r(t)+A(t,T)&+\frac{1}{2}\mid
S(t,T)\mid^2\Big)dt+ S(t,T)dW(t)\\[2ex]
&+\int_{\mathbb{R}}\Big(e^{G(t,x,T)}-1\Big)N(dt,dx)\bigg).
\end{align}
As a consequence of \eqref{rownanie na P} and definition of
$\hat{P}$ we obtain
\begin{align*}
d\hat{P}(t,T)=\hat{P}(t-,T)\bigg(\Big( A(t,T)&+\frac{1}{2}\mid
S(t,T)\mid^2\Big)dt+ S(t,T)dW(t)\\[2ex]
&+\int_{\mathbb{R}}\Big(e^{G(t,x,T)}-1\Big)N(dt,dx)\bigg).
\end{align*}
As in the case of stock market we are interested in the existence of
a martingale measure for the discounted prices. A measure $Q$ is a
martingale measure if the process $\hat{P}(\cdot,T)$ is a local
martingale with respect to $Q$ for each $T\in[0,T^\ast]$. The set of
all martingale measures is denoted by $\mathcal{Q}$. The set
$\mathcal{Q}$ is not empty if the model satisfies the $HJM$-type
conditions, that is if coefficients in  \eqref{rownanie na f} are
related in a special way. For more details see Theorem 3.13 in
\cite{BjoKabRun}. Throughout all the paper we assume that the
objective measure $P$ is at the same time a martingale one. This
assumption allows us to write the following equation for $\hat{P}$,
see Proposition 3.14 in \cite{BjoKabRun}:
\begin{align}\label{rownanie portfela}
d\hat{P}(t,T)=\hat{P}(t-,T)\bigg(S(t,T)dW(t)+\int_{\mathbb{R}}(e^{G(t,x,T)}-1)\tilde{N}(dt,dx)\bigg)
.
\end{align}

\noindent Now, let us fix a set $J$ which is assumed to be a dense,
countable subset of $[0,T^\ast]$, which elements are denoted by
$\{T_i: i\in\mathbb{N}\}$. We assume that only bonds with maturities
in $J$ are traded, i.e. only they can be used for the portfolio
construction. At the beginning we should give a precise portfolio
definition. Below it is shown a motivation for the form of the
portfolio processes used in the sequel.

Notice that if we fix $t$ then $P(t,\cdot)$, given by
\eqref{rownanie def na P}, is a continuous function on $[0,T^\ast]$,
so restricted to $J$ it is a bounded sequence. The space
\begin{gather*}
B=\Big\{z=(z_1,z_2,...): \sup_{i}\mid z_i\mid<\infty\Big\}
\end{gather*}
with the norm $\|z\|_{B}=\sup_{i}\mid z_i\mid$ is thus the state
space for the bond prices. In the classical case of stock markets
with the price process in $\mathbb{R}^d$, where $d<\infty$, it is
clear that the space of portfolios can be identified with the dual
space $(\mathbb{R}^d)^\ast=\mathbb{R}^d$. This approach is being
generalized in the context of bond markets with infinite dimensional
price process. For example in \cite{BjoKabRun} and
\cite{BjoKabRunMas} the price process takes values in
$C_0[0,\infty)$ - the space of continuous functions converging to
zero in infinity. The space of portfolios is thus
$C^\ast_0[0,\infty)$ - a space of measures with finite total
variation. In our model treating $B^\ast$ as a state space for
portfolios does not seem to be justified. The reason is that the
dual space is to large and contains abstract elements with a
doubtful financial interpretation, for example generalized Banach
limits. The portfolio space should be chosen in such a way to be
closer to practical aspects of trading. In practice the trader's
position at any time $t$ is based on finite number of bonds only, so
it is of the form
\begin{gather*}
\varphi(t)=(\varphi(t,T_{i_1}),\varphi(t,T_{i_2}),...,\varphi(t,T_{i_n}));\qquad
T_{i_j}\in J, \ j=1,2,...,n; \ n\in\mathbb{N}.
\end{gather*}
Since the number of bonds $n$ held by a trader can be arbitrarily
large, we also allow the portfolio to contain infinite number of
bonds but such that the value of the investment is finite. Since the
bond prices are bounded it is thus natural to assume that the
portfolio satisfies
\begin{gather*}
\varphi(t)=\{\varphi(t,T_{j})\}_{j=1}^{\infty};\quad
\sum_{j=1}^{\infty}\mid\varphi(t,T_{j})\mid<\infty,\qquad T_{j}\in
J, \ j=1,2,... .
\end{gather*}
Concluding, we choose $l^1\subset B^\ast$ as the portfolio space.
The value of the investment is a value of the functional
$\varphi(t)$ on the element $P(t)\in B$ and is denoted by
\begin{gather*}
<\varphi(t),P(t)>_{B^\ast,B}:=\sum_{j=1}^{\infty}\varphi(t,T_j)P(t,T_j).
\end{gather*}

\noindent By trading strategy we mean any predictable process
$\{\varphi(t); \ t\in[0,T^\ast]\}$ taking values in $l^1$. Besides
investing in bonds one can also save money in a savings account. The
wealth process at time $t$ is thus given by
\begin{gather}\label{wzor na X}
X(t)=b(t)\cdot B(t)+<\varphi(t),P(t)>_{B^{\ast},B}\qquad
t\in[0,T^\ast],
\end{gather}
where $b(t)$, $\varphi(t)$ correspond to money saved in a bank and
invested in bonds respectively. We stress the fact that dependence
on maturities $\{T_j\}_{j=1}^{\infty}$ on the right hand side of
\eqref{wzor na X} is omitted because $\varphi(t)$ and $P(t)$ are
treated as elements of infinite dimensional spaces: $l^1$ and $B$,
respectively. This notational convention will also be used with
respect to other processes appearing in the sequel.

As usual, the wealth process should be self-financing, so the
additional requirement is supposed to hold
\begin{gather}\label{wzor na dX}
dX(t)=b(t)dB(t)+<\varphi(t),dP(t)>_{B^\ast,B}\qquad t\in[0,T^\ast].
\end{gather}
Condition \eqref{wzor na dX} can be reformulated in terms of the
discounted portfolio's value. To this end we need a precise
definition of the integral
$\int<\varphi(t),d\hat{P}(t)>_{B^\ast,B}$. The definition below is
based on the equation \eqref{rownanie portfela}.
\begin{dfn}
A process $\varphi$ taking values in $l^1$ is $\hat{P}$-integrable
if it is predictable and satisfies the following conditions
\begin{gather}\label{war na calkowalnosc}
<\varphi(s),\hat{P}(s-)S(s)>_{B^\ast,B}\in\Phi, \qquad
<\varphi(s),\hat{P}(s-)(e^{G(s,x)}-1)>_{B^\ast,B}\in\Psi_{1,2}.
\end{gather}
If \eqref{war na calkowalnosc} holds, we set:
\begin{align}\label{rownanie na zdysk calke}
\int_{0}^{t}<\varphi(s),d\hat{P}(s)>_{B^\ast,B}&:=\int_{0}^{t}<\varphi(s),\hat{P}(s-)S(s)>_{B^\ast,B}dW(s)\\\nonumber
&+\int_{0}^{t}\int_{\mathbb{R}}<\varphi(s),\hat{P}(s-)(e^{G(s,x)}-1)>_{B^\ast,B}\tilde{N}(ds,dx);\quad
t\in[0,T^\ast].
\end{align}
\end{dfn}

\noindent Let us notice that integrands on the right hand side of
\eqref{rownanie na zdysk calke} are well defined since
$\hat{P}(s-)=\hat{P}(s-,\cdot)$ is a continuous function on
$[0,T^\ast]$. Indeed, due to \eqref{rownanie portfela} we obtain
$\Delta \hat{P}(t,T)=\hat{P}(t-,T)(e^{G(t,\Delta Z(t),T)}-1)$ and
putting this value to the equality
$\hat{P}(t,T)=\hat{P}(t-,T)+\Delta\hat{P}(t,T)$ we obtain
$\hat{P}(t-,T)=\frac{\hat{P}(t,T)}{e^{G(t,\Delta Z(t),T)}}$. The
last function is continuous with respect to $T$. As a consequence,
we have
\begin{gather*}
\hat{P}(t-)S(t)\in B,\qquad \hat{P}(t-)(e^{G(t,x)}-1)\in B, \qquad
\forall t\in[0,T^\ast], \quad \forall x\in\mathbb{R}.
\end{gather*}

Now, let us reformulate the self-financing condition \eqref{wzor na
dX}. Application of the integration by parts formula to the process
$\hat{X}(t):=B(t)^{-1}X(t)$ and the use of \eqref{wzor na X},
\eqref{wzor na dX} yield
\begin{align*}
d\hat{X}(t)=&B(t)^{-1}\Big(b(t)dB(t)+<\varphi(t),dP(t)>\Big)-\Big(b(t)B(t)+<\varphi(t),P(t)>\Big)B(t)^{-2}dB(t)\\[2ex]
&=<\varphi(t),B(t)^{-1}dP(t)-P(t)B(t)^{-2}dB(t)>\\[2ex]
&=<\varphi(t),d\hat{P}(t)>_{B^\ast,B}.
\end{align*}
Summarizing, the wealth process connected with a self financing
strategy can be identified with its discounted value through a pair
$(x,\varphi)$ s.t.
\begin{align*}
\hat{X}(t)=x&+\int_{0}^{t}<\varphi(s),d\hat{P}(s)>_{B^\ast,B}\\[2ex]
=x&+\int_{0}^{t}<\varphi(s),\hat{P}(s-)S(s)>_{B^\ast,B}dW(s)\\[2ex]
&+\int_{0}^{t}\int_{\mathbb{R}}<\varphi(s),\hat{P}(s-)(e^{G(s,x)}-1)>_{B^\ast,B}\tilde{N}(ds,dx);\qquad
t\in[0,T^\ast].
\end{align*}

\section{Completeness}\label{sec Completeness}
We start this section with a definition of admissible strategies - a
class of strategies involved in the definition of the market
completeness.
\begin{dfn}
Assume that a process $\varphi$ taking values in $l^1$ is $\hat{P}$-
integrable. Then $\varphi$ is an admissible strategy if the
(discounted) wealth process
\begin{gather*}
\int_{0}^{\cdot}<\varphi(s),d\hat{P}(s)>_{B^\ast,B}
\end{gather*}
is a martingale. The class of all admissible strategies will be
denoted by $\mathcal{A}$.
\end{dfn}

\noindent The definition of admissible strategies which imposes
martingale property on the wealth process is often considered in
literature, see for example \cite{Karatzas-Shreve}.

\begin{dfn} Let $A$ be a subset in the set of all
$\mathcal{F}_{T^{\ast}}$-measurable random variables interpreted as
a set of discounted contingent claims. The market is $A$-complete if
for each $X\in A$ there exists a strategy $\varphi\in\mathcal{A}$
which satisfies condition
\begin{gather}\label{war reprezentacyjny}
X=x+\int_{0}^{T^{\ast}}<\varphi(t),d\hat{P}(t)>_{B^{\ast},B},
\end{gather}
for some $x\in\mathbb{R}$. If there exists $X\in A$ s.t. condition
\eqref{war reprezentacyjny} does not hold, then the market is not
$A$-complete. If the random variable $X$ satisfies \eqref{war
reprezentacyjny} then we say that $X$ can be replicated.
\end{dfn}

\begin{lem}\label{lem o reprezentacji}
Let $\varphi\in\mathcal{A}$, $\phi\in\Phi$, $\psi\in\Psi_{1,2}$.
Assume that the proces
\begin{gather}\label{martyngalowosc sumy calek}
\int_{0}^{\cdot}\phi(s)dW(s)+
\int_{0}^{\cdot}\int_{\mathbb{R}}\psi(s,x)\tilde{N}(ds,dx)
\end{gather}
is a martingale. If the equality
\begin{gather}\label{rownosc portfela i zmiennej}
x+\int_{0}^{T^\ast}<\varphi(s),d\hat{P}(s)>_{B^\ast,B}=y+\int_{0}^{T^\ast}\phi(s)dW(s)+
\int_{0}^{T^\ast}\int_{\mathbb{R}}\psi(s,x)\tilde{N}(ds,dx)
\end{gather}
holds for some $x,y\in\mathbb{R}$ then $x=y$ and
\begin{align}\label{rownosc calki 1}
\phi(s)&=<\varphi(s),\hat{P}(s-)S(s)>_{B^\ast,B}, \qquad
P\otimes\lambda- a.s.,\\[2ex]\label{rownosc calki 2}
\psi(s,x)&=<\varphi(s),\hat{P}(s-)(e^{G(s,x)}-1)>_{B^\ast,B},\qquad
P\otimes\lambda\otimes\nu- a.s..
\end{align}
\end{lem}
{\bf Proof:} Taking expectations in \eqref{rownosc portfela i
zmiennej} we obtain $x=y$. The process
\begin{align*}
M_t:&=\int_{0}^{t}<\varphi(s),d\hat{P}(s)>_{B^\ast,B}-\int_{0}^{t}\phi(s)dW(s)-
\int_{0}^{t}\int_{\mathbb{R}}\psi(s,x)\tilde{N}(ds,dx)\\[2ex]
&=\int_{0}^{t}\left(<\varphi(s),\hat{P}(s-)S(s)>_{B^\ast,B}-\phi(s)\right)dW(s)\\[2ex]
&+\int_{0}^{t}\int_{\mathbb{R}}\left(<\varphi(s),\hat{P}(s-)(e^{G(s,x)}-1)>_{B^\ast,B}-\psi(s,x)\right)\tilde{N}(ds,dx)
\end{align*}
is thus a martingale equal to zero. With the use of Theorem \ref{tw
Kunity} we obtain \eqref{rownosc calki 1} and \eqref{rownosc calki 2}. \hfill$\square$\\

\noindent The fact of considering a specific class of admissible
strategies in the completeness problem is crucial in our approach.
If we are looking for a replicating strategy for a given integrable
random variable $X$ in the class $\mathcal{A}$ then we can identify
$X$ with a martingale $\{\mathbf{E}[X\mid\mathcal{F}_t]:
t\in[0,T^\ast]\}$. On the other hand, in view of the decomposition
\begin{gather}\label{rozklad wwo}
\mathbf{E}[X\mid\mathcal{F}_t]=\mathbf{E}X+\int_{0}^{t}\phi_X(s)dW(s)+\int_{0}^{t}\int_{\mathbb{R}}\psi_X(s,x)\tilde{N}(ds,dx),\qquad
t\in[0,T^\ast],
\end{gather}
and Theorem \ref{tw Kunity} this martingale is uniquely determined
by the processes $\phi_X,\psi_X$. Thus $X$ itself can be identified
with the integrands $\phi_X,\psi_X$. In virtue of Lemma \ref{lem o
reprezentacji} if there exists $\varphi_X\in\mathcal{A}$ satisfying
\eqref{rownosc calki 1} and \eqref{rownosc calki 2} with
$\phi=\phi_X,\psi=\psi_X$ then $\varphi_X$ is a replicating strategy
for $X$. As a consequence, if \eqref{rownosc calki 1} and
\eqref{rownosc calki 2} are not satisfied for any
$\varphi\in\mathcal{A}$ then $X$ can not be replicated.

\begin{rem}
If we do not impose any restrictions on the class of strategies or
only forbid the wealth process to take negative values then $X$ can
not be uniquely identified with the integrands $\phi_X,\psi_X$ given
by \eqref{rozklad wwo}. An example of two different integrands such
that after integrating with respect to the Wiener process give the
same bounded random variable can be found in \cite{BarJakZab},
Ex.3.10.
\end{rem}

\noindent Our method of examining conditions \eqref{rownosc calki
1}, \eqref{rownosc calki 2} is based on the following lemma which is
an extension of the moment problem solution, see Yosida \cite{Yos}.
\begin{lem}\label{lemat o momentach}
Let $\bf{E}$ be a normed linear space and $\bf{U}$ an arbitrary set.
Let $g:\bf{U}\longrightarrow \mathbb{R}$ and
$h:\bf{U}\longrightarrow \bf{E}$. Then there exists $e^\ast\in
E^\ast$ such that
\begin{gather}\label{warunek na funkcjonal}
g(u)=<e^\ast,h(u)>_{E^\ast,E}, \quad \forall u\in\bf{U},
\end{gather}
if and only if
\begin{gather}\nonumber
\exists \ \gamma>0 \quad \forall \ n\in\mathbb{N} \quad \forall \
\{\beta_i\}_{i=1}^{n}, \ \beta_i\in\mathbb{R} \quad \forall \
\{u_i\}_{i=1}^{n}, \ u_i\in{\bf U}  \quad  holds:\\[2ex]\label{warunek z momentow}
\Big|\sum_{i=1}^{n}\beta_i g(u_i)\Big|\leq\gamma \
\Big\|\sum_{i=1}^{n}\beta_i h(u_i){\Big\|}_{E}.
\end{gather}
\end{lem}
{\bf Proof:} Necessity is obvious, \eqref{warunek z momentow} holds
with $\gamma=\|e^\ast\|_{E^\ast}$. To prove sufficiency let us
define a linear subspace $\bf M$ of ${\bf E}$ by
\begin{gather*}
{\bf M}=\Big\{e\in{\bf E} : e=\sum_{i=1}^{n}\beta_ih(u_i); \quad
n\in\mathbb{N}, \ \beta_{i}\in\mathbb{R}, \ u_{i}\in{\bf U} \Big\}
\end{gather*}
and a linear transformation $\tilde{e}^\ast:{\bf M}\longrightarrow
\mathbb{R}$ by the formula
\begin{gather*}
\tilde{e}^\ast\Big(\sum_{i=1}^{n}\beta_ih(u_i)\Big)=\sum_{i=1}^{n}\beta_{i}g(u_i).
\end{gather*}
Notice, that for $e_1=\sum_{i=1}^{n}\beta_i h(u_i)$ and
$e_2=\sum_{j=1}^{m}\beta^{'}_j h(u_j)$ by \eqref{warunek z momentow}
we obtain
\begin{align*}
\Big|\tilde{e}^\ast(e_1)-\tilde{e}^\ast(e_2)\Big|&=\Big|\sum_{i=1}^{n}\beta_i
g(u_i)-\sum_{j=1}^{m}\beta^{'}_j g(u_j)\Big|\\
&\leq \gamma \Big\|\sum_{i=1}^{n}\beta_i h(u_i)
-\sum_{j=1}^{m}\beta^{'}_j h(u_j) \Big\|_{E}=\gamma \|e_1-e_2\|.
\end{align*}
If $e_1=e_2$ then $\tilde{e}^\ast(e_1)=\tilde{e}^\ast(e_2)$, so this
transformation is well defined, because its value does not depend on
the representation. It is also continuous and thus by the
Hahn-Banach theorem it can be extended to the functional $e^\ast\in
E^\ast$ which clearly satisfies \eqref{warunek na
funkcjonal}.\hfill$\square$\\

\noindent In the sequel we use the following proposition which
simplifies examining conditions \eqref{rownosc calki 1} and
\eqref{rownosc calki 2}.

\begin{prop}\label{prop o repr. produktowej}
Let $(E_1,\mathcal{E}_1,\mu_1)$, $(E_2,\mathcal{E}_2,\mu_2)$ be
measurable spaces with sigma-finite measures $\mu_1, \mu_2$ and
$(E_1\times
E_2,\mathcal{E}_1\otimes\mathcal{E}_2,\mu_1\otimes\mu_2)$ be their
product space. If two measurable functions $f_1:E_1\times
E_2\longrightarrow \mathbb{R}$, $f_2:E_1\times E_2\longrightarrow
\mathbb{R}$ satisfy condition
\begin{gather}\label{rownosc funkcji}
f_1=f_2, \qquad \mu_1\otimes\mu_2- \text{a.s.},
\end{gather}
then there exists a set $\hat{E}_1\in\mathcal{E}_1$ such that
\begin{gather}\label{1war na E1}
\hat{E}_1 \quad {\text is \ of \ full \ \mu_1 \
measure}\\\label{2war na E1} \forall x\in\hat{E}_1  \quad  {\text
the \ set}\quad \{y: f_1(x,y)=f_2(x,y)\}\quad {\text is \ of \ full
\ \mu_2 \ measure}.
\end{gather}
\end{prop}
{\bf Proof:} The assertion follows from the Fubini theorem applied
to the function $h=\mathbf{1}_{A}$ where $A:=\{(x,y)\in E_1\times
E_2 : f_1(x,y)\neq f_2(x,y)\}$.

\hfill$\square$

\subsection{L\'evy measure with a finite support}\label{subsec Levy measure with finite
support} In this section we assume that the support of the L\'evy
measure consists of finite number of points: $x_1, x_2,...,x_n$.

\noindent We start with an auxiliary lemma on linear independence of
infinite sequences. For the convenience of the reader we provide its
proof.
\begin{lem}\label{lem o liniowej niezaleznosci}
Let $M$ be an infinite matrix of the form
\begin{displaymath}
M= \left(
\begin{array}{c}
z^1\\[1ex]
z^2\\[1ex]
\vdots\\[1ex]
z^n
\end{array}
\right) = \left[
\begin{array}{cccc}
z^1_1&z^1_2&z^1_3&\ldots\\[1ex]
z^2_1&z^2_2&z^2_3&\ldots\\[1ex]
\vdots&\vdots&\vdots&\vdots\\[1ex]
z^n_1&z^n_2&z^n_3&\ldots\\[1ex]
\end{array}
\right],
\end{displaymath}
with linearly independent rows $z^1,z^2,...,z^n$. Then there exists
a set of $n$ linearly independent columns of the matrix $M$.
\end{lem}
{\bf Proof:} We will show that for some natural number $m$ the
following finite vectors
\begin{gather*}
z^k(m):=z^{k}_{1},z^{k}_{2},...,z^{k}_{m}; \qquad k=1,2,...,n,
\end{gather*}
are linearly independent. Assume, to the contrary, that for each $m$
there exist numbers $\alpha^1(m)$, $\alpha^2(m)$ ,..., $\alpha^n(m)$
such that $\sum_{k=1}^{n}\mid\alpha^{k}(m)\mid>0$ and
\begin{gather}\label{war na liniowa niezaleznosc}
\sum_{k=1}^{n}\alpha^k(m)z^k(m)=0.
\end{gather}
Without a loss of generality we can assume that
\begin{gather*}
\sum_{k=1}^{n}\mid\alpha^{k}(m)\mid=1, \qquad \forall m=1,2,... \ .
\end{gather*}
Then there exists a subsequence $m_l\rightarrow\infty$ such that
\begin{gather*}
\alpha^{k}(m_l)\longrightarrow\bar{\alpha}^k, \qquad k=1,2,...,n,
\end{gather*}
and $\sum_{k=1}^{n}\mid\bar{\alpha}^k\mid=1$. From \eqref{war na
liniowa niezaleznosc}, for each $l$, we have
\begin{gather*}
\sum_{k=1}^{n}\alpha^k(m_l)z^k(m_l)=0.
\end{gather*} Thus, for each
$\bar{m}\leq m_l$,
\begin{gather*}
\sum_{k=1}^{n}\alpha^k(m_l)z^k(\bar{m})=0.
\end{gather*}
Consequently
\begin{gather*}
\sum_{k=1}^{n}\bar{\alpha}^k z^k(\bar{m})=0, \qquad \forall
\bar{m}=1,2,... \ .
\end{gather*}
Therefore we arrive at a contradiction. \hfill $\square$

\begin{tw}\label{tw o zupelnosci nowe}
Let us assume that the following vectors in the space $B$:
\begin{gather}\label{wektory}
S(t), \ e^{G(t,x_1)}-1, \ e^{G(t,x_2)}-1, \ ... \ ,e^{G(t,x_n)}-1,
\end{gather}
are linearly independent $P\otimes\lambda$-a.s.. Then the market is
$L^1$-complete. Moreover, for each $X\in L^1$ there exists a
replicating strategy such that at any time it consists of $n+1$
bonds with different maturities.
\end{tw}
{\bf Proof:} In virtue of Lemma \ref{lem o liniowej niezaleznosci}
one can find maturities $T_{i_1}, T_{i_2},..., T_{i_{n+1}}\in J$
such that vectors
\begin{gather}
\left(
\begin{array}{c}
S(t,T_{i_j})\\[1ex]
e^{G(t,x_1,T_{i_j})}-1\\[1ex]
\vdots\\[1ex]
e^{G(t,x_n,T_{i_j})}-1
\end{array}
\right), j=1,2,...,n+1;
\end{gather}
form a set of linearly independent vectors in $\mathbb{R}^{n+1}$.
Consider any $X\in L^1$ and the representation of the process
$\mathbf{E}[X\mid\mathcal{F}_t]$ given by Theorem \ref{tw Kunity}
\begin{gather}\label{rep dla finite support}
\mathbf{E}[X\mid\mathcal{F}_t]=\mathbf{E}X+\int_{0}^{t}\phi_X(s)dW(s)+\int_{0}^{t}\int_{\mathbb{R}}\psi_X(s,x)\tilde{N}(ds,dx).
\end{gather}
\noindent Let us define a strategy $\varphi_X(t,T_{i_j});
j=1,2,...,n+1$ involving only bonds with maturities $T_{i_1},
T_{i_2},..., T_{i_{n+1}}$  as a solution of the following system of
linear equations
\begin{gather}\label{rownanie replikacyje}
\left[
\begin{array}{ccc}
S(t,T_{i_1})&...&S(t,T_{i_{n+1}})\\
e^{G(t,x_1,T_{i_1})}-1&...&e^{G(t,x_1,T_{i_{n+1}})}-1\\
\vdots&&\vdots\\
 e^{G(t,x_n,T_{i_1})}-1&...&e^{G(t,x_n,T_{i_{n+1}})}-1
\end{array}
\right] \left[
\begin{array}{c}
\hat{P}(t-,T_{i_1})\cdot{\varphi}_X(t,T_{i_1})\\
\hat{P}(t-,T_{i_2})\cdot{\varphi}_X(t,T_{i_2})\\
\vdots\\
\hat{P}(t-,T_{i_{n+1}})\cdot{\varphi}_X(t,T_{i_{n+1}})
\end{array}
\right] = \left[
\begin{array}{c}
\phi_X(t)\\
\psi_X(t,x_1)\\
\vdots\\
\psi_X(t,x_n)
\end{array}
\right]
\end{gather}
The strategy is well defined because the matrix above is
nonsingular. Moreover, ${\varphi}_X$ is a replicating strategy for
$X$. Indeed, we have
\begin{align*}
X=\mathbf{E}X&+\int_{0}^{T^\ast}\sum_{j=1}^{n+1}\hat{P}(t-,T_{i_j})S(t,T_{i_j}){\varphi}_X(t,T_{i_j})dW(t)\\&
+\int_{0}^{T^\ast}\int_{\mathbb{R}}\sum_{j=1}^{n+1}\hat{P}(t-,T_{i_j})(e^{G(t,x,T_{i_j})}-1){\varphi}_X(t,T_{i_j})\tilde{N}(dt,dx)\\
=\mathbf{E}X&+\int_{0}^{T^\ast}<{\varphi}_X(t),d\hat{P}(t)>_{B^\ast,B}.
\end{align*}
\hfill$\square$\\
\begin{rem}
It follows from the proof of Theorem \ref{tw o zupelnosci nowe} that
although the replicating portfolio contains $n+1$ different bonds,
they can change with time and are dependent on $\omega$.
\end{rem}

\begin{rem}
Theorem \ref{tw o zupelnosci nowe} shows that the assumptions of
Theorem 5.6. in \cite{BjoKabRun} can be weakened. Indeed, due to
Lemma \ref{lem o liniowej niezaleznosci} the problem is reduced to
the system of linear equations with nonsingular matrix. Thus
additional assumption imposed on coefficients $\sigma(t,\cdot)$,
$\gamma(t,\cdot)$ to be analytic functions can be relaxed. It can
also be shown that vectors \eqref{wektory} in $B$ are linearly
independent if and only if the functions: $S(t,T),e^{G(t,x_i,T)}-1,
i=1,2,...,n; T\in[0,T^\ast]$ are linearly independent. This follows
from the fact that a set of continuous functions on the interval is
linearly independent if and only if the set of their restrictions to
a fixed dense, countable subset is linearly independent. Thus
assumptions on linear independence in Theorem \ref{tw o zupelnosci
nowe} and in Theorem 5.6 in \cite{BjoKabRun} are equivalent.
\end{rem}

From practical point of view it is important to answer the question
when contingent claims can be replicated with the use of finite
number of bonds with maturities fixed at time zero. Let us notice
that if we fix $T_{i_1},T_{i_2},...,T_{i_{n+1}}\in J$ then the
method in the proof of Theorem \ref{tw o zupelnosci nowe} does not
work. The reason is that the columns of the matrix in
\eqref{rownanie replikacyje} become zero vectors if $T_{i_j}<t$ and
thus the system may not have a solution. However, we have the
following result.

\begin{tw}
Let us assume that $\sigma(t,T)$, $\gamma(t,x_i,T)$ are
deterministic functions and that the functions of variable $T$:
\begin{gather}\label{funkcje}
S(t,T), \ e^{G(t,x_1,T)}-1, \ e^{G(t,x_2,T)}-1, \ ... \
,e^{G(t,x_n,T)}-1,
\end{gather}
restricted to the set $[\bar{T},T^\ast]\cap J$, where
$\bar{T}<T^\ast$, form a set of linearly independent sequences for
almost all $t\in[0,\bar{T}]$. Moreover, assume that functions
\begin{gather*}
S(\cdot,T), \ e^{G(\cdot,x_i,T)}-1, \qquad i=1,2,...,n; \
T\in[\bar{T},T^\ast]\cap J,
\end{gather*}
are analytic on the interval $[0,\bar{T}]$. Then there exists a set
of dates $T_{i_1},T_{i_2},...,T_{i_{n+1}}\in[\bar{T},T^\ast]\cap J$,
such that each integrable, $\mathcal{F}_{\bar{T}}$ measurable random
variable $X$ can be replicated with the use of bonds with maturities
$T_{i_1},T_{i_2},...,T_{i_{n+1}}$.
\end{tw}
{\bf Proof:} Fix any $t\in[0,\bar{T}]$ such that the functions
\eqref{funkcje} restricted to $[\bar{T},T^\ast]\cap J$ are linearly
independent. In virtue of Lemma \ref{lem o liniowej niezaleznosci}
we can find maturities $T_{i_1}, T_{i_2},..., T_{i_{n+1}}\in
[\bar{T},T^\ast]\cap J$ such that the matrix
\begin{displaymath}
A(t):=\left[
\begin{array}{ccc}
S(t,T_{i_1})&...&S(t,T_{i_{n+1}})\\
e^{G(t,x_1,T_{i_1})}-1&...&e^{G(t,x_1,T_{i_{n+1}})}-1\\
\vdots&&\vdots\\
 e^{G(t,x_n,T_{i_1})}-1&...&e^{G(t,x_n,T_{i_{n+1}})}-1
\end{array}
\right]
\end{displaymath}
is invertible. Moreover, the function $det A(t), t\in[0,\bar{T}]$ is
analytic and thus can be equal to zero in a finite number of points
only. As a consequence, the matrix $A(t)$ is invertible for almost
all $t\in[0,\bar{T}]$ and thus the system
\eqref{rownanie replikacyje} has a solution on the interval $[0,\bar{T}]$.\hfill$\square$\\

\subsection{L\'evy measure with a concentration point}\label{subsec Levy measure with a concentration point}
We start examining the completeness problem in a more general
setting by introducing the following property of the L\'evy measure.
\begin{dfn}\label{def concetr. point}
The point $x_0\in\mathbb{R}$ is a concentration point of the measure
$\nu$ if there exists a sequence $\{\varepsilon_n\}_{n=1}^{\infty}$
s.t. $\varepsilon_n\searrow 0$ satisfying
\begin{gather}\label{war na punkt koncetracji}
\nu\Big\{B(x_0,\varepsilon_n)\backslash
B(x_0,\varepsilon_{n+1})\Big\}>0 \quad \forall \ n=1,2,...,
\end{gather}
where $B(x_0,\varepsilon)=\{x\in\mathbb{R}: \mid
x-x_0\mid\leq\varepsilon\}$.
\end{dfn}
Let us notice that the condition formulated in Definition \ref{def
concetr. point} is very often satisfied. For example, every L\'evy
measure with a density has a concentration point. Thus the following
theorem covers a large class of models.
\begin{tw}\label{tw_o_niezupelnosci}
Assume that the L\'evy measure $\nu$ has a concentration point
$x_0\neq 0$. If $\gamma(t,\cdot,T)$ is differentiable for each
$t\in[0,T^\ast]$, $T\in[0,T^\ast]$ and the following condition is
satisfied
\begin{gather}\label{war na rozniczkowalnosc gamma}
\qquad \forall t\in[0,T^\ast]\quad\exists \delta=\delta(t)>0\quad
s.t.\quad\int_{t}^{T^\ast}\sup_{x\in
B(x_0,\delta)}\mid\gamma_{x}^{'}(t,x,s)\mid ds<\infty
\end{gather}
then the bond market is not $L^{\infty}$-complete.
\end{tw}
{\bf Proof:} We will construct a bounded random variable $X$ which
can not be represented in the form \eqref{war reprezentacyjny} for
any strategy $\varphi\in\mathcal{A}$. At the beginning we construct
an auxiliary function $\psi$ such that
there is no $\hat{P}$-integrable process $\varphi$ satisfying condition \eqref{rownosc calki 2}.\\
Let $\{\varepsilon_n\}_{n=1}^{\infty}$ be a sequence satisfying
\eqref{war na punkt koncetracji} and define a deterministic function
$\psi$ by the formula
\[
\psi(x)=
\begin{cases}
\mid x\mid\wedge \ 1 \quad &\text{for} \quad
x\in\{B(x_0,\varepsilon_{2k+1})\backslash
B(x_0,\varepsilon_{2k+2})\} \quad k=0,1,...,
\\[2ex]
-(\mid x\mid\wedge \ 1) \quad &\text{for} \quad
x\in\{B(x_0,\varepsilon_{2k})\backslash B(x_0,\varepsilon_{2k+1})\}
\quad k=1,2,...,\\[2ex]
\mid x\mid\wedge \ 1 \quad &\text{for} \quad
x\in(-\infty,x_0-\varepsilon_1)\cup(x_0+\varepsilon_1)\cup\{x_0\}.
\end{cases}
\]

\noindent We will show that condition \eqref{rownosc calki 2} is not
satisfied by any $\hat{P}$-integrable process $\varphi$. Let us fix
any pair $(\omega,t)\in\Omega\times[0,T^\ast]$ and assume that
equality
\begin{gather}\label{war rownosc calki w dowodzie}
<\varphi(t),\hat{P}(t-)(e^{G(t,x)}-1)>_{B^\ast,B}=\psi(x)
\end{gather}
holds $\nu$ a.s.. Thus there exists a set $A_{\nu}(\omega,t)$ of a
full $\nu$ measure s.t. equality \eqref{war rownosc calki w
dowodzie} is satisfied for each $x\in A_{\nu}(\omega,t)$. Due to
Lemma \ref{lemat o momentach} there exists
$\gamma=\gamma(\omega,t)>0$ such that
\begin{gather}\nonumber
\forall \ n\in\mathbb{N} \quad \forall \ \{\beta_i\}_{i=1}^{n}, \
\beta_i\in\mathbb{R} \quad \forall \ \{x_i\}_{i=1}^{n}, \ x_i\in
A_{\nu}(\omega,t)\\[2ex]\label{war_tw}
\Big|\sum_{i=1}^{n}\beta_i\psi(x_i)\Big|\leq\gamma\Big\|
\sum_{i=1}^{n}\beta_i\hat{P}(t-)(e^{G(t,x_i)}-1)\Big\|_{B}.
\end{gather}
Let us notice that due to \eqref{war na punkt koncetracji} we have
\begin{gather*}
\nu\Big\{A_{\nu}(\omega,t)\cap \big\{B(x_0,\varepsilon_n)\backslash
B(x_0,\varepsilon_{n+1})\big\}\Big\}>0
\end{gather*}
so we can choose a sequence $\{a_k\}_{k=1}^{\infty}$ s.t.
\begin{gather*}
a_k\in A_{\nu}(\omega,t)\cap \big\{B(x_0,\varepsilon_k)\backslash
B(x_0,\varepsilon_{k+1})\big\} \quad  \forall \ k=1,2,... .
\end{gather*}
\noindent Let us examine the condition \eqref{war_tw} with $n=2$,
$\beta_1=1, \beta_2=-1$ and $x_1=a_{2k+1}$, $x_2=a_{2k+2}$ for
$k=0,1,...$. Then the left hand side of \eqref{war_tw} is of the
form
\begin{gather*}
\frac{1}{\gamma}\Big|\beta_1\psi(a_{2k+1})+\beta_2\psi(a_{2k+2})\Big|=\frac{1}{\gamma}
\Big((\mid a_{2k+1}\mid \wedge \ 1)+(\mid a_{2k+2} \mid \wedge \
1)\Big)
\end{gather*}
and thus satisfies
\begin{gather*}
\lim_{k\longrightarrow\infty}\frac{1}{\gamma}\Big|\beta_1\psi(a_{2k+1})+\beta_2\psi(a_{2k+2})\Big|=\frac{2(\mid
x_0\mid\wedge \ 1)}{\gamma}\neq 0.
\end{gather*}
\noindent In estimating of the right hand side of \eqref{war_tw} we
will use the inequality \eqref{oszacowanie G} and \eqref{oszacowanie G prim} below.\\
In view of \eqref{war na rozniczkowalnosc gamma} we have
\begin{align}\label{oszacowanie G}\nonumber
\sup_{T\in J}\sup_{x\in B(x_0,\delta)}&\mid G(t,x,T)\mid
\leq\sup_{T\in J}\sup_{x\in
B(x_0,\delta)}\int_{t}^{T}\mid\gamma(t,x,s)\mid ds\\\nonumber
&\leq\sup_{T\in J}\int_{t}^{T}\sup_{x\in
B(x_0,\delta)}\mid\gamma(t,x,s)\mid ds\\\nonumber &\leq\sup_{T\in
J}\int_{t}^{T}\Big\{\mid\gamma(t,x_0,s)\mid+\sup_{x\in
B(x_0,\delta)}\mid \gamma_{x}^{'}(t,x,s)\mid2\delta \Big\}ds\\
&\leq\int_{t}^{T^\ast}\mid\gamma(t,x_0,s)\mid
ds+2\delta\int_{t}^{T^\ast}\sup_{x\in B(x_0,\delta)}\mid
\gamma_{x}^{'}(t,x,s)\mid ds<\infty.
\end{align}
The condition \eqref{war na rozniczkowalnosc gamma} implies
differentiability of $G(t,\cdot,T)$ and the following estimation
\begin{align}\label{oszacowanie G prim}\nonumber
\sup_{T\in J}\sup_{x\in B(x_0,\delta)}\mid
G_{x}^{'}(t,x,T)\mid&=\sup_{T\in J}\sup_{x\in
B(x_0,\delta)}\mid \int_{t}^{T}\gamma_{x}^{'}(t,x,s)ds\mid\\
&\leq\int_{t}^{T^\ast}\sup_{x\in
B(x_0,\delta)}\mid\gamma_{x}^{'}(t,x,s)\mid ds<\infty.
\end{align}

\noindent The right hand side of \eqref{war_tw} can be estimated as
follows
\begin{align*}
\Big\|\hat{P}(t-)(e^{G(t,a_{2k+1})}-1)&-\hat{P}(t-)(e^{G(t,a_{2k+2})}-1)\Big\|_{B}\\[2ex]
&= \sup_{T\in
J}\Big|\hat{P}(t-,T)(e^{G(t,a_{2k+1},T)}-1)-\hat{P}(t-,T)(e^{G(t,a_{2k+2},T)}-1)\Big|\\[2ex]
&\leq \sup_{T\in J}|\hat{P}(t-,T)| \ \sup_{T\in J}
\Big|e^{G(t,a_{2k+1},T)}-e^{G(t,a_{2k+2},T)}\Big|.
\end{align*}
The first supremum is finite since $\hat{P}(t-,\cdot)$ is a
continuous function. To deal with the second supremum let us notice
that for sufficiently large $k$ the points $a_{2k+1},a_{2k+2}$ are
in $B(x_0,\delta)$ and thus we have
\begin{align}\label{nier na pochodna}\nonumber
\sup_{T\in J}
\Big|e^{G(t,a_{2k+1},T)}-e^{G(t,a_{2k+2},T)}\Big|\leq\sup_{T\in J
}\sup_{x\in
B(x_0,\delta)}\Big|\frac{d}{dx}e^{G(t,x,T)}\Big|\cdot\mid
a_{2k+1}-a_{2k+2}&\mid\\[2ex]
\leq\sup_{T\in J}\sup_{x\in B(x_0,\delta)}e^{\mid G(t,x,T)\mid}
\cdot \sup_{T\in J}\sup_{x\in B(x_0,\delta)}\mid
G_{x}^{'}(t,x,T)\mid \cdot\mid a_{2k+1}-a_{2k+2}&\mid .
\end{align}
In view of \eqref{oszacowanie G} and \eqref{oszacowanie G prim} we
see that the last factor in \eqref{nier na pochodna} goes to $0$
when
$k\rightarrow\infty$.\\

\noindent Thus we conclude that condition \eqref{war_tw} is not
satisfied for any $(\omega,t)\in \Omega\times[0,T^\ast]$ and thus
\eqref{war rownosc calki w dowodzie} does not hold $\nu - a.s.$ for
any $(\omega,t)\in \Omega\times[0,T^\ast]$. As a consequence of
Proposition \ref{prop o repr. produktowej} there is no
$\hat{P}$-integrable process satisfying \eqref{rownosc calki 2}.\\

\noindent Now, with the use of the function $\psi$, we construct a
bounded random variable $X$ which can not be replicated. \\

\noindent It is clear that  $\psi\in\Psi_{1,2}$. Let us define the
stopping time $\tau_k$ by
\begin{gather*}
\tau_k=\inf\{t:
\Big|\int_{0}^{t}\int_{\mathbb{R}}\psi(x)\tilde{N}(ds,dx)\Big|\geq
k\}\wedge T^{\ast}
\end{gather*}
and choose a number $k_0$ s.t. the set
$\{(\omega,\tau_{k_0}(\omega));
\omega\in\Omega\}\subseteq\Omega\times[0,T^\ast]$ is of positive
$P\otimes\lambda$ measure. Then the process
$\psi(x)\mathbf{1}_{(0,\tau_{k_0}]}(s)$ is predictable and bounded.
The random variable
\begin{gather}\label{zmienna niereplikowalna ograniczona}
X=\int_{0}^{T^{\ast}}\int_{\mathbb{R}}\psi(x)\mathbf{1}_{(0,\tau_{k_0}]}(s)\tilde{N}(ds,dx)
\end{gather}
is thus well defined and it is also bounded because
$|\Delta\int_{0}^{\cdot}\int_{\mathbb{R}}\psi(x)\tilde{N}(ds,dx)|\leq
1$. For any $(\omega,t)\in \{(\omega,\tau_{k_0}(\omega));
\omega\in\Omega\}$ condition \eqref{war_tw} is not satisfied
$\nu$-a.s.. As a consequence of Proposition \ref{prop o repr.
produktowej} condition \eqref{rownosc calki 2} is not satisfied by
any $\hat{P}$- integrable process. Moreover,
$\int_{0}^{\cdot}\int_{\mathbb{R}}\psi(s,x)\tilde{N}(ds,dx)$ is a
martingale. As a consequence of Lemma \ref{lem o reprezentacji}
there is no admissible strategy which replicates $X$.\hfill$\square$

\subsection{L\'evy measure with a discrete support}\label{subsec Levy measure with discrete
support} In this section we consider the L\'evy measure with a
support consisting of infinite number of discrete points denoted by
$\{x_i\}_{i=1}^{\infty}$. To exclude the case studied in Section
\ref{subsec Levy measure with a concentration point} we assume that
the support has no concentration point, so the sequence satisfies
\begin{gather}\label{war na nosnik w disrete support}
\lim_{i\rightarrow\infty}\mid x_i\mid=\infty.
\end{gather}
\noindent Let us notice, that in this case the L\'evy measure is a
sequence of positive numbers $\{\nu(x_i)\}_{i=1}^{\infty}$ which,
due to relation $\int_{\mathbb{R}}\mid x\mid\wedge \ 1 \
\nu(dx)<\infty$, satisfies condition
\begin{gather}\label{war na miare Levyego dla discrete support}
\sum_{i=1}^{\infty}\nu(\{x_i\})<\infty.
\end{gather}

\noindent In the following theorem we show that under additional
condition imposed on the coefficient $\gamma$ we obtain a result on
incompleteness.
\begin{tw}\label{tw o niezup w discrete suppport 1}
Assume that the following set
\begin{gather*}
A=\Big\{(\omega,t)\in\Omega\times[0,T^{\ast}] \quad s.t. \quad
G(t,x_i,T)\leq0 \quad \forall T\in[0,T^\ast]\quad \forall i=1,2,...
\Big\}
\end{gather*}
is of positive $P\otimes\lambda$ measure. Then the market is not
$L^2$-complete.
\end{tw}
{\bf Proof:} We construct a random variable $X\in L^2$ which can not
be represented in the form \eqref{war reprezentacyjny}. At the
beginning, using condition \eqref{war na miare Levyego dla discrete
support}, let us define a sequence $\{\psi(x_i)\}_{i=1}^{\infty}$
which depends neither on $\omega$ nor $t$ in the following way
\begin{equation}\label{wzor na psi rosnace do infty}
\psi(x_i)=
\begin{cases}
\sqrt{k} &\text{for} \quad i=i_k
\\
0 &\text{for} \quad i\neq i_k,
\end{cases}
\end{equation}
\noindent where $i_k:=\inf\Big\{i: \nu(x_i)\leq\frac{1}{k^3}\Big\}$.
This sequence satisfies the following two conditions
\begin{align}\label{war pierwszy na ciag}
\lim\sup_{i\rightarrow\infty}\mid \psi(x_i)\mid=\infty,\\
\label{war drugi na ciag}
\sum_{i=1}^{\infty}\mid\psi(x_i)\mid^2\nu(\{x_i\})\leq\sum_{k=1}^{\infty}\frac{1}{k^2}<\infty.
\end{align}
We show that the representation \eqref{rownosc calki 2} which we
write in the form
\begin{gather}\label{war reprezentacyjny dla discrete support}
<\varphi(t),\hat{P}(t-)(e^{G(t,x_i)}-1)>_{B^\ast,B}=\psi(x_i) \quad
\forall i=1,2,...,
\end{gather}
does not hold $P\otimes\lambda\otimes\nu$-a.s. for any $\hat{P}$
integrable process $\varphi$. Let us fix $(\omega,t)\in A$ and
assume to the contrary that \eqref{war reprezentacyjny dla discrete
support} is satisfied for some $\varphi(t)$. Then by Lemma
\ref{lemat o momentach} there exists $\gamma=\gamma(\omega,t)>0$
such that
\begin{gather}\nonumber
\forall \ n\in\mathbb{N} \quad \forall \ \{\beta_{k}\}_{k=1}^{n}, \
\beta_{k}\in\mathbb{R} \quad \forall \ \{x_{i_k}\}_{k=1}^{n}\\[2ex]\label{war rownowazny dla discrete support}
\Big|\sum_{k=1}^{n}\beta_{k}\psi(x_{i_k})\Big|\leq\gamma\Big\|
\sum_{k=1}^{n}\beta_{k}\hat{P}(t-)(e^{G(t,x_{i_k})}-1)\Big\|_{B}.
\end{gather}
Let us check \eqref{war rownowazny dla discrete support} with $n=1$,
$\beta_1=1$ and for $i_1=1,2,...$ successively, that is
\begin{gather}\label{war uproszczony dla discrete support}
\Big|\psi(x_{i})\Big|\leq\gamma\sup_{T\in J}\Big|
\hat{P}(t-,T)(e^{G(t,x_{i},T)}-1)\Big| \qquad \forall i=1,2,... .
\end{gather}
By the definition of the set $A$ for any $i=1,2,...$ we have
\begin{gather*}
\mid e^{G(t,x_{i},T)}-1\mid\leq1 \quad \forall \ T\in J.
\end{gather*}
Using the inequality
\begin{gather*}
\sup_{T\in J}\Big|
\hat{P}(t-,T)(e^{G(t,x_{i},T)}-1)\Big|\leq\sup_{T\in
J}\Big|\hat{P}(t,T)\Big|\cdot\sup_{T\in
J}\Big|e^{G(t,x_{i},T)}-1\Big|
\end{gather*}
and the fact that $\hat{P}(t,\cdot)$ is continuous we see that
\begin{gather*}
\underset{i\rightarrow\infty}{\lim\sup}\sup_{T\in J}\Big|
\hat{P}(t-,T)(e^{G(t,x_{i},T)}-1)\Big|<\infty.
\end{gather*}
However, recall that  the left hand side of \eqref{war uproszczony
dla discrete support} satisfies \eqref{war pierwszy na ciag}, so the
required constant $\gamma$ does not exist. We have shown that for
any $(\omega,t)\in A$ the representation \eqref{war reprezentacyjny
dla discrete support} does not hold. But $P\otimes\lambda(A)>0$, so
in view of Proposition \ref{prop o repr. produktowej}, the
representation \eqref{war reprezentacyjny dla discrete support} does
not hold
$P\otimes\lambda\otimes\nu$-a.s. for any $\hat{P}$-integrable process.\\

\noindent In view of \eqref{war drugi na ciag} we see that
$\psi\in\Psi_{1,2}$ and that the process
$\int_{0}^{\cdot}\int_{\mathbb{R}}\psi(x)\tilde{N}(ds,dx)$ is a
martingale. Thus with the use of Lemma \ref{lem o reprezentacji} we
conclude that the following random variable
\begin{gather}\label{niereplikowalna z L2 w discrete support}
X:=\int_{0}^{T^\ast}\int_{\mathbb{R}}\psi(x)\tilde{N}(ds,dx)
\end{gather}
can not be replicated by strategies from the class $\mathcal{A}$. By
application isometric formula to $X$ we obtain that $X$ is square
integrable.
\hfill$\square$\\

\noindent The next theorems are based on the behavior of the
expression $\parallel G(t,x_i)\parallel_B$ for large $i$. Since
their proofs are similar to those presented earlier, we provide the
sketches only.
\begin{tw}
If the following condition holds
\begin{gather}\label{war na granice 0}
\underset{\mid x_i\mid\rightarrow\infty}{\lim\inf}
\parallel G(t,x_i)\parallel_B=0, \qquad P\otimes\lambda-a.s.
\end{gather}
then the market in not $L^\infty$-complete.
\end{tw}
{\bf Proof:} The condition \eqref{war na granice 0} implies
\begin{gather*}
\underset{\mid x_i\mid\rightarrow\infty}{\lim\inf}\parallel
e^{G(t,x_i)}-1\parallel_B\leq \underset{\mid
x_i\mid\rightarrow\infty}{\lim}e^{\parallel G(t,x_i)\parallel_B}-1=0
\end{gather*}
For $\psi(x_i)\equiv1$ condition \eqref{war uproszczony dla discrete
support} is thus not satisfied what we can check by calculating
$\lim\inf _{i}$ for both sides.

\noindent The bounded random variable which can not be replicated is
constructed in the same way as in the proof of Theorem
\eqref{tw_o_niezupelnosci}, see formula \eqref{zmienna
niereplikowalna ograniczona}. \hfill$\square$

\begin{tw}
If the set
\begin{align}\label{war na granice alpha}\nonumber
A=\Big\{(\omega,t)\in\Omega\times[0,T^\ast]&: \exists
\alpha=\alpha(\omega,t); 0<\alpha<\infty \ s.t.\\[2ex]
&\lim_{\mid x_i\mid\rightarrow\infty}\parallel
G(t,x_i)\parallel_B=\alpha \Big\}
\end{align}
is of positive $P\otimes\lambda$ measure then the market in not
$L^2$-complete.
\end{tw}
{\bf Proof:} We use $\psi$ constructed in the proof of Theorem
\ref{tw o niezup w discrete suppport 1}, given by the formula
\eqref{wzor na psi rosnace do infty}. Then \eqref{war na granice
alpha} implies that
\begin{gather*}
\underset{\mid x_i\mid\rightarrow\infty}{\lim\sup} \
\frac{\mid\psi(x_i)\mid}{\parallel e^{G(t,x_i)}-1\parallel_B}=\infty
\end{gather*}
and thus condition \eqref{war uproszczony dla discrete support} does
not hold. A square integrable random variable which can not be
replicated is given by \eqref{niereplikowalna z L2 w discrete
support}. \hfill$\square$\\

\noindent To study the case when $\parallel G(t,x_i)\parallel_B$
tends to infinity we restrict ourselves to the linear form of the
coefficient $\gamma$, i.e. $\gamma(t,x,T)=\gamma(t,T)x$. This is
done to simplify a formulation of the next theorem. Notice that in
this case we have $G(t,x,T)=G(t,T)x$.

\begin{tw} Assume that $\gamma(t,x,T)=\gamma(t,T)x$.
If there exists a constant $\tilde{G}>0$ such that the set
\begin{gather*}
A=\Big\{(\omega,t)\in\Omega\times[0,T^\ast]: \ \parallel
G(t,T)\parallel_B\leq \tilde{G} \Big\}
\end{gather*}
is of positive $P\otimes\lambda$ measure and the L\'evy measure has
exponential moment of order $2(\tilde{G}+\varepsilon)$ for some
$\varepsilon>0$, i.e.
\begin{gather*}
\sum_{i=1}^{\infty}e^{2(\tilde{G}+\varepsilon)\mid
x_i\mid}\nu(\{x_i\})<\infty,
\end{gather*}
then the market is not $L^2$-complete.
\end{tw}
{\bf Proof:} Define
\begin{gather*}
\psi(x_i)=e^{(\tilde{G}+\varepsilon)\mid x_i\mid} ,\qquad i=1,2,...
.
\end{gather*}
For any $(\omega,t)\in A$, condition \eqref{war uproszczony dla
discrete support} is not satisfied because we have
\begin{gather*}
\lim_{i\rightarrow\infty}\frac{\mid\psi(x_i)\mid}{\parallel e^{G(t)
x_i}-1\parallel_B}\geq
\lim_{i\rightarrow\infty}\frac{\mid\psi(x_i)\mid}{\mid
e^{\tilde{G}\mid x_i\mid}-1\mid}=\infty.
\end{gather*}

\noindent As a consequence the following random variable
\begin{gather*}
X:=\int_{0}^{T^\ast}\int_{\mathbb{R}}\psi(x)\tilde{N}(ds,dx)
\end{gather*}
can not be replicated and it is square integrable because
\begin{gather*}
\mathbf{E}(X^2)=\mathbf{E}\int_{0}^{T^\ast}\sum_{i=1}^{\infty}e^{2(\tilde{G}+\varepsilon)\mid
x_i\mid}\nu(\{x_i\})ds<\infty.
\end{gather*}
\hfill$\square$

\begin{rem}
In this paper we assume that only bonds with maturities in $J$ can
be traded and thus we accepted $B$ for the state space. However, if
we admit for the portfolio construction all bonds with maturities in
$[0,T^\ast]$ and the state space $C([0,T^\ast])$ - a space of
continuous functions with the supremum norm, then all the results
remain true. This is because for any continuous function
$h:[0,T^\ast]\longrightarrow \mathbb{R}$ we have
\begin{gather*}
\|h\|_B=\sup_{T\in J}\mid h(T)\mid=\sup_{T\in [0,T^\ast]}\mid
h(T)\mid=\|h\|_{C([0,T^\ast])}
\end{gather*}
and thus all the arguments based on the norm in $B$ can be
automatically replaced by the norm in $C([0,T^\ast])$.
\end{rem}


\begin{thebibliography}{9}

\bibitem{Apl}
Applebaum, D.: "L\'evy Processes and Stochastic Calculus", (2004),
Cambridge University Press,


\bibitem{BarJakZab}
Barski M., Jakubowski J., Zabczyk J.: "On incompleteness of bond
markets with infinite number of random factors", (2011),
Mathematical Finance, 21,3, 541-556,


\bibitem{BjoKabRunMas}
Bj\"ork T., Di Masi G., Kabanov Y. and Runggaldier W.: "Towards a
general theory of bond markets", (1997), {\it Finance and
Stochastics} 1, 141-174,

\bibitem{BjoKabRun}
Bj\"ork T., Kabanov Yu., Runggaldier W.: "Bond market structure in
the presence of marked point process", (1997), {\it Mathematical
Finance}, 7, 211-239,

\bibitem{Car}
Carmona R. and Tehranchi M.: "A characterization of hedging
portfolios for interest rate contingent claims", (2004),  {\it
Annals of Applied Probability}, 14, 1267-1294,

\bibitem{EbJacRaib}
Eberlein E., Jacod J., Raible S.: "L\'evy term structure models:
No-arbitrage and completeness", (2005), {\it Finance and
Stochastics}, 9, 67-88,


\bibitem{Karatzas-Shreve}
Karatzas, I., Shreve, S.E.: "Methods of mathematical finance",
(1998), Springer,


\bibitem{Kun}
Kunita H.: "Representation of martingales with jumps and
applications to mathematical finance", (2004), {\it Advanced Studies
in Pure Mathematics} 41, {\it Stochastic Analysis and Related
Topics}, 209-232,

\bibitem{Taf}
Taflin E.: "Bond market completeness and attainable contingent
claims", (2005), {\it Finance and Stochastics} 9, 429 - 452,

\bibitem{Yos}
Yosida K.: "Functional Analysis", (1980), Springer.






\end{thebibliography}
\end{document}